\documentclass[11pt,reqno]{amsart}

\theoremstyle{plain}
\newtheorem{theorem} {Theorem}[section]
\newtheorem{lemma}[theorem] {Lemma}
\newtheorem{proposition}[theorem] {Proposition}
\newtheorem{corollary}[theorem] {Corollary}
\newtheorem{property}[theorem] {Property}

\theoremstyle{definition}

\theoremstyle{remark}
\newtheorem{remark}[theorem] {Remark}

\numberwithin{equation}{section}

\usepackage{amsfonts}
\usepackage{amssymb}
\usepackage{amscd}

\newcommand{\R}{{\mathbb R}}
\newcommand{\Z}{{\mathbb Z}}
\newcommand{\N}{{\mathbb N}}

\newcommand{\PP}{{\mathcal P}}

\newcommand{\CC}{{\mathbb C}}

\newcommand{\GG}{{\mathfrak G}}
\newcommand{\al}{{\alpha}}

\newcommand{\sa}{{\sigma}}

\newcommand{\iy}{{\infty}}
\newcommand{\vphi}{{\varphi}}
\newcommand{\vep}{{\varepsilon}}
\newcommand{\g}{{\gamma}}
\newcommand{\de}{{\delta}}

\newcommand{\z}{{\zeta}}
\newcommand{\be}{{\beta}}
\newcommand{\G}{{\Gamma}}

\newcommand{\bna}{\begin{eqnarray}}
\newcommand{\ena}{\end{eqnarray}}
\newcommand{\ba}{\begin{eqnarray*}}
\newcommand{\ea}{\end{eqnarray*}}
\newcommand{\beq}{\begin{equation}}
\newcommand{\eeq}{\end{equation}}

\begin{document}

\title[Asymptotic Relations]
{Asymptotic Relations Between interpolation Differences
and Zeta Functions}
\author{Michael I. Ganzburg}
 \address{212 Woodburn Drive, Hampton,
 VA 23664\\USA}
 \email{michael.ganzburg@gmail.com}

 \begin{abstract}
 Asymptotic relations
between zeta functions
(such as, $\zeta(s),\,\beta(s)$,
 and  other Dirichlet $L$-functions)
and interpolation differences
 of functions like $\vert y\vert^s$ and their interpolating
entire functions of exponential type $1$ are discussed.
New criteria for zeros of the zeta
 functions
  in the critical strip in terms of integrability of
the interpolation differences are obtained as well.
 \end{abstract}
 \maketitle

 \section{Introduction}\label{S1}
\setcounter{equation}{0}
In this paper we find asymptotic relations
between zeta functions
(such as, $\zeta(s),\,\beta(s)$,
 and  other Dirichlet $L$-functions)
and interpolation differences
 of functions like $\vert y\vert^s$ and their interpolating
entire functions of exponential type $1$.
As corollaries, we  obtain new criteria for zeros of the zeta
 functions
  in the critical strip in terms of integrability of
the interpolation differences.

\subsection{Notation}\label{SS1.1}
Let
$\Z$ denote the set of all integers
and let $\Z_+$ be the set
of all nonnegative integers
and $\N:=\Z_+\setminus\{0\}$.
Let $K(R):=\{w\in\CC:\vert w\vert= R\}$
be the circle in $\CC=\R+i\R$ centered at the origin
of radius $R>0$.
In addition, $\vert \Omega\vert$ denotes the Lebesgue measure
of a  measurable set $\Omega\subset\R$.
We also use the floor function
 $\lfloor a \rfloor,\,a\ge 0$,
 and the gamma function $\G(s),\,s\in\CC$.

Let $L_p(\Omega)$ be the space of all measurable
complex-valued functions $F$
 on a measurable set $\Omega\subseteq\R$  with the finite quasinorm
 \ba
 \|F\|_{L_p(\Omega)}:=\left\{\begin{array}{ll}
 \left(\int_\Omega\vert F(y)\vert^p\, dy\right)^{1/p}, & 0<p<\iy,\\
 \mbox{ess} \sup_{y\in \Omega} \vert F(y)\vert, &p=\iy.
 \end{array}\right.
 \ea

 We say that an entire function $g:\CC\to \CC$ has
 exponential type $\sa$
 if for any $\vep>0$ there exists a constant
 $C=C(\vep,g,\sa)>0$ such that
 for all $z\in \CC$,
 $\vert g(z)\vert\le C\exp\left((\sa+\vep)\vert z\vert\right)$.
  The class of all entire functions of exponential type $\sa$
  is denoted by $B_\sa$.
  Throughout the paper, if no confusion may occur,
  the same notation is applied to
  $g\in B_\sa$ and its restriction to $\R$ (e.g., in the form
  $g\in  B_\sa\cap L_p(\R))$.
  Here, we mostly discuss
  entire functions of exponential type $1$ (EFET1) from $B_1$.

  The Fourier transform of a function $G\in L_1(\R)$
   is denoted by the formula
   \ba
   \widehat{G}(u):=
   \int_{\R}G(y)e^{iuy}dy.
   \ea
\noindent
We use the same notation $\widehat{G}$
for the Fourier transform of a tempered distribution
 $G$ on $\R$.
 By the definition (see, e.g., \cite[Sect. 4.1]{S2003}),
 $G$ is a continuous linear functional $\langle G,\psi\rangle$ on the Schwartz
 class $S(\R)$ of all
 test functions $\psi$ on $\R$, and $\widehat{G}$ is defined by the formula
 $\langle \widehat{G},\psi\rangle
 :=\langle\psi,\widehat{\psi}\rangle,\,\psi\in S(\R).$
 If a function $G$ of polynomial growth on $\R$ is locally
 integrable, then it generates the tempered
 distribution $G$ by the formula
 $ \langle G,\psi \rangle:=\int_\R G(y)\psi(y)dy,\,
 \psi\in S(\R)$.

Throughout the paper $C,\,C_1,\,C_2,\ldots,C_{18}$
denote positive constants independent
of essential parameters.
 Occasionally, we indicate dependence on certain parameters.
 The same symbol $C$ does not
 necessarily denote the same constant in different occurrences,
 while $C_j,\,1\le j\le 18$,
 denotes the same constant in different occurrences.

 \subsection{Interpolation Differences
and Zeta Functions}\label{SS1.2}
It was Bernstein  who in 1938 initiated the study of
polynomial approximation
and approximation by EFET1 of the function
$f_s(y):=\vert y\vert^s$
 by proving the following celebrated result
(see \cite[Eqn. (36)]{B1938}) for $s>0$:
\beq\label{E1.1}
\lim_{n\to\iy}n^s \inf_{P_n\in\PP_{n}}
\left\|
f_s-P_n\right\|_{L_\iy([-1,1])}
=\inf_{g\in B_1}
\left\|
f_s-g\right\|_{L_\iy(\R)}
<\iy,
\eeq
where $\PP_n$ is the class of all univariate algebraic
polynomials of degree at most $n$.
An essential ingredient of the proof of
\eqref{E1.1} in \cite{B1938} was the use of
the interpolation difference
$\Delta_s(y):=f_s(y)-g_s(y),\,y\in\R,$
where $g_s$ is the only EFET1 that
interpolates $f_s$ at the nodes
$\{\pi(n+1/2)\}_{n\in\Z}$ under
a certain condition.
Bernstein \cite[Eqn. (42)]{B1938} announced without
proof the formula for $g_s$ and the integral
representation for $\Delta_s$.
The author
\cite[Lemma 5 (a)]{G2002} proved these formulae and
obtained similar results for
$\vert y\vert^s\mathrm{sgn}\,y$.
Note that $L_p$-versions of \eqref{E1.1}
 for complex $s$ with $\mathrm{Re}\,s>
 \max\{-1,-1/p\},\,p\in(0,\iy]$, or $\mathrm{Re}\,s=0,\,p=\iy$,
  were recently proved in \cite{G2021a}.

  The systematic studies  of the interpolation difference
  $f-g$ and its $L_1(\R)$-norm for certain real-valued functions $f$
  (including $f=f_s, \,s>0$) were conducted by
  Vaaler \cite{V1985}, Littmann \cite{L2005},
  Carneiro and Vaaler \cite{CV2010}, the author \cite{G2010},
  and others.

It turns out that the interpolation
difference $\Delta_s$
possesses the following surprising
property:
\beq\label{E1.2}
\be(s)=\frac{\pi}{4\sin(\pi s/2)}
\lim_{\vert y\vert\to\iy}
\frac{\Delta_s(y)}{\cos y},
\qquad s>0,\quad s\ne 2,\,4,\ldots,
\eeq
where $\be(s)$ is the Dirichlet beta function
(or the Dirichlet $L$-function $L(s,\chi)$
with the character $\chi$ of modulus $4$, see
Section 1.3 for definitions).
This relation can be extended to a complex $s$
as well.

In this paper we discuss more general interpolation
differences
$\Delta_{k,s,v}(y):=f_{k,s,v}(y)-g_{k,s,v}(y),\,y\in\R,$
related to the function $\Phi$
(see
Section 1.3 for the definition),
with a parameter $v\in(0,1]$ and a complex $s,\,
\mathrm{Re}\,s>1$ for $k=0$ and
$\mathrm{Re}\,s>0$ for $k=1$.
Here, $f_{k,s,v}$ is the linear combination of
$\vert y\vert^s$ and
$\vert y\vert^s\mathrm{sgn}\,y$,
and $g_{k,s,v}$ is the only EFET1 that interpolates
$f_{k,s,v}$ at the nodes
$\{\pi(n+k/2)\}_{\vert n\vert=1-k}^\iy$
under certain conditions.

Main results are presented in Section \ref{S2}.
The explicit formulae for
$f_{k,s,v}$ and $g_{k,s,v}$
and a general version
of \eqref{E1.2} are given in
Theorem \ref{T2.1}.
Special cases associated with the zeta functions
$\zeta,\,\beta$ and general
$L$-functions are discussed
in Corollary \ref{C2.3} and Theorem \ref{T2.7},
respectively.
Different versions of Theorem \ref{T2.7}
for two special $L$-functions
are presented in Corollaries \ref{C2.10}
and \ref{C2.9}.
New criteria for zeros of the zeta
 functions and general $L$-functions in the critical
 strip in terms of integrability of
the interpolation differences are discussed in
Corollaries \ref{C2.2}, \ref{C2.4}, \ref{C2.5}, and
\ref{C2.8}.

 Preliminaries are discussed below,
 and the proofs of main results are given in Section \ref{S5}.
 The proofs are based on the two lemmas proved in
 Sections \ref{S3} and \ref{S4}.

 \subsection{Preliminaries}\label{SS1.3}
 Here, we discuss certain special functions and their properties.
 \vspace{.12in}\\
 \textbf{Special Functions $\Phi,\zeta$, and $\be$.}
 The Lerch transcendent
 \ba
 \Phi(z,s,v):=\sum_{n=0}^\iy(v+n)^{-s}z^n,
 \quad\vert z\vert<1,\quad v\ne 0,\,-1,\ldots,\quad s\in\CC,
 \ea
 of three complex variables can be extended to a different domain
 by the following integral representation
 \beq\label{E2.1}
 \Phi(z,s,v)=\frac{1}{\G(s)}\int_0^\iy\frac{t^{s-1}e^{(1-v)t}}
 {e^t-z}dt,
 \eeq
 where $\mbox{Re}\,v>0$ and either
 $\vert z\vert\le 1,\,z\ne 1,\, \mbox{Re}\,s>0$, or
 $z=1,\,\mbox{Re}\,s>1$
 (see \cite[Sect. 1.11]{Erd1953}, \cite{LG2002}, \cite{LL2012}).

 In this paper we discuss $\Phi(z,s,v)$ for $v>0$ (mostly for
 $v\in(0,1]$) and
 either
 $z=1,\, \mbox{Re}\,s>1$, or $z=-1,\,\mbox{Re}\,s>0$,
 because zeta functions can be expressed in terms of this function.
 The Hurwitz (or generalized) zeta function
 $\z(s,v):=\sum_{n=0}^\iy(v+n)^{-s},\,\mbox{Re}\,s>1,\, v>0$;
 the Riemann zeta function
 $\z(s):=\sum_{n=1}^\iy n^{-s},\,\mbox{Re}\,s>1$;
 and the Dirichlet beta function
 $\be(s):=\sum_{n=0}^\iy (-1)^n (2n+1)^{-s},\,\mbox{Re}\,s>0$,
 allow the following representations in terms of $\Phi(\pm 1,s,v)$
 (see \cite[Sects. 1.10 and 1.12]{Erd1953}):
 \bna
&&\z(s,v)= \Phi(1,s,v),\quad \mbox{Re}\,s>1,
\quad v>0;\label{E2.2}\\
&&\z(s)= \Phi(1,s,1), \quad \mbox{Re}\,s>1;
\label{E2.3}\\
&&\z(s)= \left(2^s-1\right)^{-1}\Phi(1,s,1/2),
\quad \mbox{Re}\,s>1;\label{E2.5}\\
&&\z(s)= \left(1-2^{1-s}\right)^{-1}\Phi(-1,s,1),
\quad \mbox{Re}\,s>0,
\quad s\ne 1;\label{E2.4}\\
&&\be(s)= 2^{-s}\Phi(-1,s,1/2), \quad \mbox{Re}\,s>0.
\label{E2.6}
 \ena
 Note that formula \eqref{E2.4} extends $\z(s)$ to
 $\mbox{Re}\,s>0,\,
 s\ne 1$.\vspace{.12in}\\
 \textbf{The Dirichlet Characters and $L$-functions.}
 Let $\chi=\chi(\cdot,q):\Z\to\CC$ be a \emph{Dirichlet
character of modulus $q$} with $q\in\N,\,q>1$.
Then
(see, e.g., \cite[Sect. 4.2]{P1957})
$\chi$ is a completely multiplicative and
$q$-periodic function on $\Z$
with $\vert\chi(\cdot)\vert=0$ or $1$,
and $\chi(l)=0$
if and only if $(l,q)>1$. In addition,
\beq\label{E2.6.1}
\sum_{l=1}^{q-1}\chi(l)=0,\qquad \chi\ne \chi_0,
\eeq
where
\ba
\chi_0(l)=\chi_0(l,q):=\left\{\begin{array}{ll}
1, &(l,q)=1,\\
0, &(l,q)>1, \end{array}\right.\qquad l\in\Z,
\ea
is the \emph{principal character}.
Two more examples are given below.
\beq\label{E2.6.1a}
\chi(l,3)=\left\{\begin{array}{rlr}
1,&l\equiv 1\, (\mathrm{mod}\,3),\\
-1,&l\equiv 2 \,(\mathrm{mod}\,3),\\
0,&l\equiv 0 \,(\mathrm{mod}\,3),
\end{array}\right.\,\,
\chi(l,4)=\left\{\begin{array}{ll}
(-1)^{(l-1)/2}, &l \,\,\mathrm{is\,\,odd},\\
0, &l \,\,\mathrm{is\,\,even},
\end{array}\right.\,\, \chi\ne\chi_0,\,\,
l\in\Z.
\eeq
A Dirichlet character $\chi$ of modulus $q$  is
called \emph{primitive} if for every proper divisor $d$
of $q$ (that is, $d<q$), there exists an integer
$a\equiv 1\,(\mathrm{mod}\,d)$,
with $(a,q)= 1$ and $\chi(a)\ne 1$.
The following properties of primitive characters
hold true (see, e.g., \cite[Sect. 6.3]{C2012}):
\begin{itemize}
\item[(i)]
If $q$ is
an odd prime, then every nonprincipal character
is primitive.
\item[(ii)]
If $\chi$ is a primitive character,
then
\ba
\left\vert\sum_{l=1}^{q-1}\chi(l)e^{2\pi i l/q}
\right\vert=\sqrt{q}.
\ea
 \item[(iii)]
If $\chi$ is a primitive character,
then
 \ba
 \sum_{l=1}^{q-1}\chi(l)e^{2\pi i nl/q}
 =\Bar{\chi}(n)
 \sum_{l=1}^{q-1}\chi(l)e^{2\pi i l/q},
 \qquad n\in\Z,
 \ea
 where $\Bar{\chi}(n)$ is the complex conjugate of
 $\chi(n)$.
 Note that the sum in the right-hand
 side of this equality
  is called the \emph{Gauss sum}.
\end{itemize}
 Let us consider the following equation in
 $n\in\Z$:
\beq\label{E2.6.2}
\sum_{l=1}^{q-1}\chi(l)e^{2\pi i nl/q}=0.
\eeq

\begin{proposition}\label{P1.1}
  Let $\chi\ne\chi_0$ be a primitive character (for example,
  $\chi \ne \chi_0$ and $q$ is an
odd prime by property
(i)). Then the following statements hold true:\\
(a) A number $n\in\Z\setminus\{0\}$
satisfies equation \eqref{E2.6.2}
if and only if $(n,q)>1$.\\
(b) The function
\beq\label{E2.6.2a}
T_q(2z):=\frac{\sin q z}
{\sum_{l=1}^{q-1}\chi(l)e^{i(2l-q) z}},
\qquad z\in\CC,\quad q>2,
\eeq
is entire if and only if $q=3$ or $q=4$.
In addition,
\beq\label{E2.6.2b}
T_3(y)=(i/2)(1+2\cos y),\qquad T_4(y)=(i/2)\cos y.
\eeq
\end{proposition}
\proof
Statement (a) immediately follows from properties
(ii) and (iii) of primitive characters.
To prove statement (b), we note first that
$T_3$ and $T_4$, given by
\eqref{E2.6.2b}, are entire functions. Next,
the exponential sum
\beq\label{E2.6.2c}
P_{q}(2z):=(2/i)\sum_{l=1}^{q-1}\chi(l)e^{i(2l-q) z}
=(2/i)e^{-q z}\sum_{l=1}^{q-1}\chi(l)e^{2li z}
\eeq
is a trigonometric polynomial of exact degree $q-2$
(since $\chi(1)\ne 0$), and
$\left\vert P_{q}(2\cdot)\right\vert$ is
a $\pi$-periodic function on $\CC$ by \eqref{E2.6.2c}.
 Then $P_{q}(2\cdot)$  has exactly $q-2$ zeroes
in the strip $\mathrm{Re}\,z\in(0,\pi]$.
Furthermore, by statement (a), the number of zeros of
$P_{q-2}$ of the form $\pi n/q,\,n\in\N,\,n\le q$,
is $q-\vphi(q)$, where $\vphi$
is  Euler's totient function. Therefore, if
$T_q(2\,\cdot)$ is entire, then $\vphi(q)=2$.
This is possible only for $q\in\{3,\,4,\,6\}$.
Finally, $\chi(\cdot,3)$ and $\chi(\cdot,4)$
 are primitive characters and
 $\chi(\cdot,6)$ is an imprimitive character.
 Thus statement (b) is established.
 \hfill $\Box$\vspace{.12in}\\
 Note that all above-mentioned Dirichlet characters
 have modulus $q>1$. However, the Dirichlet character
 of modulus $1$ can be
 defined as well by $\chi=\chi(\cdot,1)\equiv 1$
  on $\Z$.

A Dirichlet $L$-function $L(s,\chi)$ is a meromorphic
function on $\CC$, which is the
holomorphic extension of
a Dirichlet $L$-series
 \beq\label{E2.6.3}
L(s,\chi):=\sum_{l=1}^\iy\frac{\chi(l)}{l^s},
\qquad \mathrm{Re}\,s>1.
\eeq
In particular (see, e.g., \cite[Sect. 4.3]{P1957}),
\beq\label{E2.6.3a}
L(s,\chi(\cdot,1))=\zeta(s),\quad
L(s,\chi(\cdot,2))=L(s,\chi_0(\cdot,2)), \quad
L(s,\chi_0(\cdot,q))
=\z(s)\prod_{p\vert q}\left(1-p^{-s}\right),
\eeq
where
$s\in\CC\setminus \{1\},\,q>1$, and $p$ is a prime.
Two more examples below follow from \eqref{E2.6.1a}
and \eqref{E2.6.3}:
\beq\label{E2.6.3b}
L(s,\chi(\cdot,3))=\sum_{d=0}^\iy
\left((3d+1)^{-s}-(3d+2)^{-s}\right),\,\,
L(s,\chi(\cdot,4))
=\be(s),
\,\,\, \mathrm{Re}\,s>0,\,\,\chi\ne\chi_0.
\eeq
In general, if $\chi\ne\chi_0$, then $L(s,\chi)$ is
an entire function and series  \eqref{E2.6.3}
is convergent for $\mathrm{Re}\,s>0$.
In addition, the following integral
representation for $\mathrm{Re}\,s>0$ holds true:
\beq\label{E2.6.4}
L(s,\chi)=\frac{1}{\G(s)}\int_0^\iy\frac{x^{s-1}}
{e^{qx}-1} \sum_{l=1}^{q-1}\chi(l)e^{(q-l) x}dx,
\qquad \chi\ne \chi_0,\quad \mathrm{Re}\,s>0.
\eeq
This formula easily follows from \eqref{E2.6.3} for
$\mathrm{Re}\,s>1$, and its holomorphic extension to
$\mathrm{Re}\,s>0$ immediately follows from
\eqref{E2.6.1}.
Representation \eqref{E2.6.4} can be rewritten in
terms of the function $\Phi$
(or in terms of $\zeta(s,v)$, see \eqref{E2.2})
by using \eqref{E2.1} as
\beq\label{E2.6.5}
L(s,\chi)=q^{-s}\sum_{l=1}^{q-1}\chi(l)\Phi(1,s,l/q),
\qquad \mathrm{Re}\,s>0.
\eeq
Note again that despite the fact that $\Phi(1,s,l/q)$
is defined for $\mathrm{Re}\,s>1$, the $L$-function
in \eqref{E2.6.5} can be holomorphically extended to
$\mathrm{Re}\,s>0$ by \eqref{E2.6.4}.

\section{Main Results}\label{S2}
\setcounter{equation}{0}
Here, we present limit representations for
$\Phi(\pm 1,s,v),\,\zeta(s),\,\be(s)$,
 and general $L$-functions and apply them
to new criteria for their zeros.

\subsection{Asymptotic Relations for
$\Phi(\pm 1,s,v)$.}\label{S2.1}
 Throughout the paper
 we assume that $k=0$ or $k=1$ and

\beq\label{E2.6d}
m_k=m_k(s):=\left\{\begin{array}{ll}
\lfloor (\mbox{Re}\,s-1)/2\rfloor,&k=0,\\
\lfloor (\mbox{Re}\,s)/2\rfloor,&k=1.
\end{array}\right.
\eeq
We first discuss the following general results about
 $\Phi\left((-1)^k,s,v\right)$ and its zeros.

\begin{theorem} \label{T2.1}
Let one of the following conditions on $k,\,s$, and $v$
be satisfied:
\bna
&&\mathrm{Re}\,s\in (1,\iy),\,\,
\mathrm{Re}\,s\ne 3,\,5,\ldots,\,\,
\mathrm{and}\,\,s\notin\N,\,\,
 \mathrm{if}\,\, k=0\,\, \mathrm{and}\,\, v\in (0,1);
 \label{E2.6a}\\
 &&\mathrm{Re}\,s\in \cup_{d\in\Z_+} (1+2d,2+2d),\,\,
 \mathrm{if}\,\, k=0\,\, \mathrm{and}\,\, v=1;
 \label{E2.6b}\\
 &&\mathrm{Re}\,s\in (0,\iy),\,\,
 \mathrm{Re}\,s\ne 2,\,4,\ldots,\,\,\mathrm{and}\,\,
 \,s\notin\N,\,\,
 \mathrm{if}\,\, k=1\,\, \mathrm{and}\,\, v\in (0,1]
  \label{E2.6c}.
 \ena
Then the following statements hold true:\\
(a) The function
\bna\label{E2.7}
g_{k,s,v}(y)&:=&-\sin(y+\pi k/2)
\left(\sum_{j=1}^{m_k}(-1)^{j-1}\G(s-2j)
\Phi\left((-1)^k,s-2j,v\right)(2y)^{2j}
\right.\nonumber\\
&+&\left.\frac{\pi 2^{s}y^{2(m_k+1)}}{\sin(\pi s)}
\sum_{n=1-k}^\iy
\frac{[\pi(n+k/2)]^{s-2m_k-1}\sin[\pi(2n+k)v+\pi s/2]}
{y^2-[\pi(n+k/2)]^2}
\right)
\ena
is the only EFET1 that interpolates the function
\beq\label{E2.8}
f_{k,s,v}(y):=-\pi 2^{s-2}
\left(\frac{\vert y\vert^s\sin[(2v-1)y-\pi k/2]}
{\sin(\pi s/2)}
+\frac{\vert y\vert^s\mathrm{sgn}\,y\cos[(2v-1)y-\pi k/2]}
{\cos(\pi s/2)}\right)
\eeq
at the nodes $\{\pi(n+k/2)\}_{n\in\Z}$
and satisfies the following conditions:
(C1) $f_{k,s,v}-g_{k,s,v}\in L_\iy(\R)$ and
(C2) $g_{k,s,v}^{(1-k)}(0)=0$.\\
(b) In addition to statement (a), for
$n\in\Z$,
\bna\label{E2.8a}
&&f_{k,s.v}[\pi(n+k/2)]=g_{k,s.v}[\pi(n+k/2)]\nonumber\\
 &&=-\frac{\pi 2^{s-1}(-1)^{n+k}[\pi
 \vert n+k/2\vert]^s}
 {\sin(\pi s)}
 \sin[\pi(2n+k)v+(\pi s/2)\,\mathrm{sgn}(n+k/2)].
 \ena
(c) The limit equality
\beq\label{E2.9}
\Phi\left((-1)^k,s,v\right)
=
\frac{1}{\G(s)}\lim_{\vert y\vert\to\iy}\frac{f_{k,s,v}(y)-g_{k,s,v}(y)}
{\sin(y+\pi k/2)}
\eeq
is valid.
\end{theorem}

\begin{corollary} \label{C2.2}
Let $k,\,s,$ and $v$ be as in Theorem \ref{T2.1}.\\
(a) If
$\Phi\left((-1)^k,s,v\right) =0$,
then for any $p\in(1/2,\iy),\,
f_{k,s,v}-g_{k,s,v}\in L_p(\R)$.\\
(b) If there exists $p\in(1/2,\iy)$ such that
$f_{k,s,v}-g_{k,s,v}\in L_p(\R)$, then
$\Phi\left((-1)^k,s,v\right) =0$.
\end{corollary}

\subsection{Asymptotic Relations for
$\z$ and $\be$.}\label{S2.2}
Here, we present simplified formulae
for $\z(s)$ and $\be(s)$
with $\mbox{Re}\,s$ belonging to
 one of the intervals $(1,2),\,(1,3),\,(0,2)$.
Note that for these values of $s$ the polynomial
in \eqref{E2.7} is
zero. Then the next corollary follows directly from
Theorem \ref{T2.1} (c) and equalities
\eqref{E2.3} through \eqref{E2.6}.

\begin{corollary} \label{C2.3}
The following statements hold true:\\
(a) For $\mathrm{Re}\,s\in(1,2)$,
\ba
\z(s)&=&\frac{1}{\G(s)}
\lim_{\vert y\vert\to\iy}\frac{f_{0,s,1}(y)-g_{0,s,1}(y)}
{\sin y}\nonumber\\
&=&\frac{\pi 2^{s-2}}{\G(s)}\lim_{\vert y\vert\to\iy}
\frac{1}{\sin y}\nonumber\\
&\times&\left(
-\frac{\vert y\vert^s\sin y}
{\sin(\pi s/2)}
-\frac{\vert y\vert^s\mathrm{sgn}\,y\,\cos y}
{\cos(\pi s/2)}
+\frac{2y^{2}\sin y}{\cos(\pi s/2)}
\sum_{n=1}^\iy
\frac{(\pi n)^{s-1}}
{y^2-(\pi n)^2}
\right).
\ea
(b) For $\mathrm{Re}\,s\in(1,3)$,
\ba
\z(s)&=&\frac{1}{\left(2^s-1\right)\G(s)}
\lim_{\vert y\vert\to\iy}
\frac{f_{0,s,1/2}(y)-g_{0,s,1/2}(y)}
{\sin y}\nonumber\\
&=&\frac{\pi 2^{s-2}}{\left(2^s-1\right)\G(s)\cos(\pi s/2)}
\lim_{\vert y\vert\to\iy}
\frac{1}{\sin y}\nonumber\\
&\times&\left(
-\vert y\vert^s\mathrm{sgn}\,y
+2y^{2}\sin y
\sum_{n=1}^\iy
\frac{(-1)^n(\pi n)^{s-1}}
{y^2-(\pi n)^2}
\right).
\ea
(c) For $\mathrm{Re}\,s\in(0,2)$,
\bna\label{E2.12}
\z(s)&=&\frac{1}{\left(1-2^{1-s}\right)\G(s)}
\lim_{\vert y\vert\to\iy}\frac{f_{1,s,1}(y)-g_{1,s,1}(y)}
{\cos y}\nonumber\\
&=&\frac{\pi 2^{s-2}}{\left(1-2^{1-s}\right)\G(s)}
\lim_{\vert y\vert\to\iy}
\frac{1}{\cos y}\nonumber\\
&\times&\left(
\frac{\vert y\vert^s\cos y}
{\sin(\pi s/2)}
-\frac{\vert y\vert^s\mathrm{sgn}\,y\,\sin y}
{\cos(\pi s/2)}
-\frac{2y^{2}\cos y}{\cos(\pi s/2)}
\sum_{n=0}^\iy
\frac{[\pi (n+1/2)]^{s-1}}
{y^2-[\pi (n+1/2)]^2}
\right).
\ena
(d) For $\mathrm{Re}\,s\in(0,2)$,
\bna\label{E2.13}
\be(s)&=&\frac{1}{2^s\G(s)}
\lim_{\vert y\vert\to\iy}
\frac{f_{1,s,1/2}(y)-g_{1,s,1/2}(y)}
{\cos y}\nonumber\\
&=&\frac{\pi}{4\G(s)\sin(\pi s/2)}
\lim_{\vert y\vert\to\iy}
\frac{1}{\cos y}\nonumber\\
&\times&\left(
\vert y\vert^s
+2y^{2}\cos y
\sum_{n=0}^\iy
\frac{(-1)^n[\pi (n+1/2)]^{s-1}}
{y^2-[\pi (n+1/2)]^2}
\right).
\ena
\end{corollary}
\noindent

\begin{remark}\label{R2.3a}
The formulae for $f_{k,s,1/2}$ and
$g_{k,s,1/2}$ from statements (b) and (d)
of Corollary \ref{C2.3} for real $s$
were obtained in \cite[Lemma 5 (a)]{G2002}.
These results for $k=1$ were announced in
\cite[Eqn. (42)]{B1938}.
\end{remark}

The specified criteria, presented
in the next two corollaries,
follow directly from
 Corollary \ref{C2.2} and equalities
 \eqref{E2.4} and \eqref{E2.6}.
Note that the formulae for
the interpolation differences
$f_{1,s,1}-g_{1,s,1}$
and $f_{1,s,1/2}-g_{1,s,1/2}$
are given in \eqref{E2.12}
 and \eqref{E2.13},
respectively.

\begin{corollary} \label{C2.4}
Let $\mathrm{Re}\, s\in(0,1).$\\
(a) If
$\z(s) =0$, then for any $p\in(1/2,\iy),\,
f_{1,s,1}-g_{1,s,1}\in L_p(\R)$.\\
(b) If there exists $p\in(1/2,\iy)$ such that
$f_{1,s,1}-g_{1,s,1}\in L_p(\R)$, then
$\z(s) =0$.
\end{corollary}

\begin{corollary} \label{C2.5}
Let $\mathrm{Re}\, s\in(0,1).$\\
(a) If
$\be(s) =0$, then for any $p\in(1/2,\iy),\,
f_{1,s,1/2}-g_{1,s,1/2}\in L_p(\R)$.\\
(b) If there exists $p\in(1/2,\iy)$ such that
$f_{1,s,1/2}-g_{1,s,1/2}\in L_p(\R)$, then
$\be(s) =0$.
\end{corollary}

\begin{remark}\label{R2.6}
There are dozens of the other criteria
(see, e.g.,\cite{C2003}).
However, the criteria
of Corollaries \ref{C2.4} and \ref{C2.5}
are the first in terms of
interpolation of functions like
$\vert y\vert^s$ by EFET1.
Note that certain criteria in terms of
interpolation of functions like
$\vert y\vert^s$ by algebraic polynomials are
presented in \cite[Sects. 5.3 and 5.4]{G2013}.
\end{remark}

\subsection{Asymptotic Relations for
Dirichlet $L$-functions.}\label{S2.3}
Recall that
$m_1=m_1(s)=\lfloor (\mbox{Re}\,s)/2\rfloor$
by \eqref{E2.6d}.
Let $\chi:\Z\to\CC$ be a Dirichlet
character of modulus $q$ with $q>1$ and
$\chi\ne\chi_0$,
and let $E(\chi)$ be a
symmetric set of all
nonzero integers $n$,
satisfying equation \eqref{E2.6.2}.
Note that by \eqref{E2.6.1}, $0\in E(\chi)$.
In particular, by Proposition \ref{P1.1} (a),
in case of a primitive character $\chi,\,
n\in E(\chi)$ if and only if
$(n,q)>1$.
In addition, let $L(s,\chi)$ be a
Dirichlet $L$-function.
Then the following results are valid.

\begin{theorem} \label{T2.7}
Let
$\mathrm{Re}\,s\in (0,\iy),\,
 \mathrm{Re}\,s\ne 2,\,4,\ldots $, and
 $s\notin\N.$
Then the following statements hold true:\\
(a) The function
\bna\label{E2.13a}
\g_{s,q}(y)
&:=&-\sin y
\left(q^s\sum_{j=1}^{m_1}(-1)^{j-1}\G(s-2j)
L\left(s-2j,\chi\right)(2y)^{2j}
\right.\nonumber\\
&+&\left.\frac{\pi 2^{s}y^{2(m_1+1)}}{\sin(\pi s)}
\sum_{n\in \N\setminus E(\chi)}
\frac{(\pi n)^{s-2m_1-1}
\sum_{l=1}^{q-1}\chi(l)\sin[2\pi n l/q+\pi s/2]}
{y^2-(\pi n)^2}
\right)
\ena
is the only EFET1 that interpolates the function
\beq\label{E2.13b}
\vphi_{s,q}(y):=-\pi 2^{s-2}
\left(\frac{\vert y\vert^s
\sum_{l=1}^{q-1}\chi(l)\sin[(2l/q-1)y]}
{\sin(\pi s/2)}
+\frac{\vert y\vert^s\mathrm{sgn}\,y
\sum_{l=1}^{q-1}\chi(l)\cos[(2l/q-1)y]}
{\cos(\pi s/2)}\right)
\eeq
at the nodes $\{\pi n\}_{n\in\Z}$
and satisfies the following conditions:
(C1*) $\vphi_{s,q}
-\g_{s,q}\in L_\iy(\R)$ and
(C2*) $\g_{s,q}^{\prime}(0)=0$.\\
(b) In addition to statement (a),
\bna\label{E2.13c}
&&\vphi_{s,q}(\pi n)
=\g_{s,q}(\pi n)\nonumber\\
 &&=\left\{\begin{array}{ll}
 0, &n\in E(\chi),\\
 -\frac{\pi 2^{s-1}(-1)^{n}[\pi
 \vert n\vert]^s}
 {\sin(\pi s)}
 \sum_{l=1}^{q-1}\chi(l)
 \sin[2\pi n l/q+(\pi s/2)\,\mathrm{sgn}\,n],
 &\vert n\vert\ge 1, \,n\in\Z\setminus E(\chi).
 \end{array}\right.
 \ena
(c) The limit equality
\beq\label{E2.13d}
L\left(s,\chi\right)
=
\frac{1}{q^s\G(s)}\lim_{\vert y\vert\to\iy}
\frac{\vphi_{s,q}(y)-\g_{s,q}(y)}
{\sin y}
\eeq
is valid.
\end{theorem}

\begin{remark}\label{R2.7a}
Though Theorem \ref{T2.7} holds true
for $q>1$ and $\chi\ne\chi_0$, certain versions of
the theorem are valid for $q=1$ or
$\chi=\chi_0$.
Indeed, \eqref{E2.3} and the
first formula of
 \eqref{E2.6.3a} show that
 if condition \eqref{E2.6b} on $s$ is satisfied,
 then Theorem \ref{T2.1}
 (see also Corollary \ref{C2.3} (a)
 for $\mathrm{Re}\,s\in(1,2)$)
 can substitute Theorem \ref{T2.7} for $q=1$.
 In addition, it follows from
 the third formula of \eqref{E2.6.3a}
 that the same conclusion can be made
for $q>1$ and $\chi=\chi_0$
(up to the constant in \eqref{E2.6.3a}).
 In particular,
 statement (b) of Corollary \ref{C2.3}
 follows from relation
 \eqref{E2.13d} for $q=2$ and $m_1=0$
 by \eqref{E2.6.3a}.
 \end{remark}

\begin{corollary} \label{C2.8}
Let $\mathrm{Re}\, s\in(0,1).$\\
(a) If
$L(s,\chi) =0$, then for any $p\in(1/2,\iy),\,
\vphi_{s,q}-\g_{s,q}\in L_p(\R)$.\\
(b) If there exists $p\in(1/2,\iy)$ such that
$\vphi_{s,q}-\g_{s,q}\in L_p(\R)$, then
$ L(s,\chi)=0$.
\end{corollary}
One of the shortcomings of Theorem \ref{T2.7}
is the same set of nodes
$\{\pi n\}_{\vert n\vert=1}^\iy$ for all $q>1$.
Assuming that $T_q\in B_1$,
we can obtain more informative results
when we replace
in Theorem \ref{T2.7}
$\vphi_{s,q}(y)$ with
 $\vphi_{s,q}^*(y):=\vphi_{s,q}(qy/2)/P_{q}(y)$
 and
$\g_{s,q}(y)$ with
 $\g_{s,q}^*(y):=\g_{s,q}(qy/2)/P_{q}(y)$,
where $T_{q}(2\,\cdot)$ and
$P_{q}(2\,\cdot)$ are defined by
\eqref{E2.6.2a} and \eqref{E2.6.2c},
respectively.
The advantage of using $\vphi_{s,q}^*$ and
$\g_{s,q}^*$ instead of $\vphi_{s,q}(y)$ and
$\g_{s,q}(y)$ is due to the fact that
the EFET1 $\g_{s,q}^*$ interpolates $\vphi_{s,q}^*$
only at the nodes
$\{2\pi n/q\}_{n\in\Z\setminus E(\chi)}$,
unlike \eqref{E2.13c}.

Note that by Proposition \ref{P1.1} (b),
$T_q\in B_1$ if and only if
$q=3$ or $q=4$ with $T_3$ and $T_4$
given in \eqref{E2.6.2b}
and
$P_3(y)=\sin(y/2),
\,P_4(y)=2\sin y$.
In the following two corollaries
we obtain more explicit results when compared with
Theorem \ref{T2.7} for $q=3$ and $q=4$.
A more explicit version of Theorem \ref{T2.7}
for $q=3$ is presented below.

\begin{corollary} \label{C2.10}
If $s$ satisfies the conditions
 of Theorem \ref{T2.7},
then the following statements hold true:\\
(a) The function
\bna\label{E2.14a}
&&\g_{s,3}^*(y):=(1+2\cos y)
\left(-3^s\sum_{j=1}^{m_1}(-1)^{j-1}\G(s-2j)
L\left(s-2j,\chi(\cdot,3)\right)(3y)^{2j}
\right.\nonumber\\
&&+\left.\frac{\pi 3^{s-1/2}y^{2(m_1+1)}}{\sin(\pi s/2)}
\left(-\sum_{d=0}^\iy
\frac{[2\pi (d+1/3)]^{s-2m_1-1}}
{y^2-[2\pi (d+1/3)]^2}
+\sum_{d=0}^\iy
\frac{[2\pi (d+2/3)]^{s-2m_1-1}}
{y^2-[2\pi (d+2/3)]^2}\right)
\right)
\ena
is the only EFET1 that interpolates the function
\beq\label{E2.14b}
\vphi_{s,3}^*(y)
:=\frac{\pi(3\vert y\vert)^s}{2\sin(\pi s/2)}
\eeq
at the nodes $\{2\pi (d\pm 1/3)\}_{d\in\Z}$
and satisfies the following conditions:
(C1*) $\vphi_{s,3}^*
-\g_{s,3}^*\in L_\iy(\R)$ and
(C2*) $\g_{s,3}^{*}(0)=0$.\\
(b) In addition to statement (a),
\beq\label{E2.14c}
\vphi_{s,3}^*[2\pi (d\pm 1/3)]
=\g_{s,3}^*[2\pi (d\pm 1/3)]
 =\frac{\pi 3^{s}[2\pi
 \vert d\pm 1/3\vert]^s}
 {2\sin(\pi s/2)},
 \qquad d\in\Z.
 \eeq
(c) The limit equality
\beq\label{E2.14d}
L\left(s,\chi(\cdot,3)\right)
=
\frac{1}{3^s\G(s)}\lim_{\vert y\vert\to\iy}
\frac{\vphi_{s,3}^*(y)-\g_{s,3}^*(y)}
{1+2\cos y}
\eeq
is valid.
\end{corollary}
The next corollary shows that for $q=4$
Theorem \ref{T2.7} can be replaced by more
 explicit Theorem \ref{T2.1} for $k=1$
 (see also Corollary \ref{C2.3} (d)
 for $\mathrm{Re}\, s\in(0,2)$).

\begin{corollary} \label{C2.9}
If $s$ satisfies the conditions
 of Theorem \ref{T2.7}, then
\beq\label{E2.14e}
\vphi_{s,4}^*=2^s f_{1,s,1/2},\qquad
\g_{s,4}^*=2^s g_{1,s,1/2},
\eeq
and
\beq\label{E2.14f}
L\left(s,\chi(\cdot,4\right))
=
\frac{1}{4^s\G(s)}\lim_{\vert y\vert\to\iy}
\frac{\vphi_{s,4}^*(y)-\g_{s,4}^*(y)}
{\cos y}.
\eeq
\end{corollary}

\begin{remark}\label{C2.11}
Versions of Corollary \ref{C2.8}
for $q=3$ and $q=4$
with
$\vphi_{s,q}$ and $\g_{s,q}$
replaced by
 $\vphi_{s,q}^*$ and $\g_{s,q}^*$,
 respectively,
 hold true. In case of $q=4$ this version
  is equivalent to Corollary \ref{C2.5}.
\end{remark}

The proofs of Theorems \ref{T2.1} and \ref{T2.7}
and Corollaries \ref{C2.2}, \ref{C2.8}, \
\ref{C2.10}, and \ref{C2.9} are presented in
Section \ref{S5}.
The proof of Theorem \ref{T2.1}
is based on the two lemmas that are proved in
Sections \ref{S3} and \ref{S4}.

\section{Properties of the Integral}\label{S3}
\setcounter{equation}{0}
Recall that $\Phi,\,f_{k,s,v},$ and
$g_{k,s,v}$ are defined by \eqref{E2.1},
\eqref{E2.8}, and \eqref{E2.7},
respectively.

\subsection{Four Major Properties.}\label{S3.1}
The proof of Theorem \ref{T2.1}
  is based
on properties of the integral
\beq\label{E2.14}
F_{k,s,v}(y):=
\int_0^\iy\frac{t^{s-1}e^{(1-v)t}}
 {\left(e^t-(-1)^k\right)
 \left(1+[t/(2y)]^2\right)}dt,
\eeq
where $y\in\R\setminus\{0\}$ is a fixed number,
$k=0$ or $k=1,\,v>0,$ and $\mathrm{Re}\,s>1$ if $k=0$
and $\mathrm{Re}\,s>0$ if $k=1$.
Note that these conditions guarantee the
absolute convergence
 of the integral and the boundedness of
 $F_{k,s,v}$ on $\R$
 (that is, $F_{k,s,v}\in L_\iy(\R)$).

\begin{lemma}\label{L2.7}
If $y\in\R\setminus\{0\}$,
then the following statements hold true.\\
(a) The following inequality is valid:
\beq\label{E2.15}
\left\vert F_{k,s,v}(y)-
\G(s)\Phi\left((-1)^k,s,v\right)\right\vert
\le C_1(k,s,v)y^{-2}.
\eeq
(b) If  $\vert y\vert\le 1,\,
\mathrm{Re}\,s\in(2,\iy)$,
and $v\in(0,1]$, then
\beq\label{E2.15a}
\vert F_{0,s,v}(y)\vert\le C_2(s,v)
\vert y\vert^{2(\mathrm{Re}\,s-1)/\mathrm{Re}\,s}.
\eeq
(c) For any $m\in \Z_+$
such that $\mathrm{Re}\,s-2m>1$ if $k=0$
and $\mathrm{Re}\,s-2m>0$ if $k=1$, the function
$F_{k,s,v}$ satisfies the recurrence relation
\beq\label{E2.16}
F_{k,s,v}(y)
=\sum_{j=1}^{m}(-1)^{j-1}\G(s-2j)
\Phi\left((-1)^k,s-2j,v\right)(2y)^{2j}
+(-1)^m(2y)^{2m}F_{k,s-2m,v}(y).
\eeq
(d) Let $k,\,s$, and $v$ satisfy one of conditions
\eqref{E2.6a}, \eqref{E2.6b}, and \eqref{E2.6c}.
 Then the series
 \beq\label{E2.16a}
 \frac{\pi 2^{s}}{\sin(\pi s)}
 \sum_{n=1-k}^\iy
[\pi(n+k/2)]^{s-2m_k-1}\sin[\pi(2n+k)v+\pi s/2]
\frac{\sin(y+\pi k/2)}
{y^2-[\pi(n+k/2)]^2}
\eeq
  in the definition
\eqref{E2.7} of $g_{k,s,v}$ is convergent
for every $y\in\R\setminus\{0\}$, and
the following representation holds true:
\beq\label{E2.17}
\sin(y+\pi k/2)F_{k,s,v}(y)
=f_{k,s,v}(y)-g_{k,s,v}(y).
\eeq
\end{lemma}

\subsection{Proofs of Statements (a),
 (b), and (c).}\label{S3.2}
(a) Inequality \eqref{E2.15} immediately follows
from the equality
\ba
\G(s)\Phi\left((-1)^k,s,v\right)
-F_{k,s,v}(y)
=\int_0^\iy\frac{t^{s+1}e^{(1-v)t}}
 {\left(e^t-(-1)^k\right)
 \left(t^2+4y^2\right)}dt.
\ea
(b) Setting $\al:=2/\mathrm{Re}\,s$, we split
\beq\label{E2.17a}
\vert F_{0,s,v}(y)\vert
\le
  \left(\int_0^{y^\al}
 +\int_{y^\al}^1
 +\int_1^\iy\right)
 \frac{t^{\mathrm{Re}\,s-1}
e^{(1-v)t}}
 {\left(e^t-1\right)
 \left(1+[t/(2y)]^2\right)}dt
 =I_1(y)+I_2(y)+I_3(y),
\eeq
where
\bna
&&I_1(y)\le
e^{1-v}\int_0^{y^\al}t^{\mathrm{Re}\,s-2}dt
\le \left(e^{1-v}/(\mathrm{Re}\,s-1)
\right)y^{\al(\mathrm{Re}\,s-1)};\label{E2.17b}\\
&&I_2(y)\le
4e^{1-v}\int_{y^\al}^1
\frac{t^{\mathrm{Re}\,s-3}}
{t/y^2}dt
\le 4\left(e^{1-v}/(\mathrm{Re}\,s-2)
\right)y^{2-\al};\label{E2.17c}\\
&&I_3(y)\le 8y^2\int_1^\iy
t^{\mathrm{Re}\,s-3}e^{-vt}dt
\le 8\G(\mathrm{Re}\,s-2)y^2.\label{E2.17d}
\ena
Combining \eqref{E2.17b}, \eqref{E2.17c}, and
 \eqref{E2.17d}
with \eqref{E2.17a}, we arrive at
\eqref{E2.15a}.
\vspace{.12in}\\
(c) Using the elementary identity ($h\in\R,\,u\in\R$)
\ba
{h^m=(h+u)\sum_{j=1}^{m}
(-1)^{j-1}h^{m-j}u^{j-1}+(-1)^mu^m}
\ea
for $h=t^2$ and $u=(2y)^2$, we obtain
\bna\label{E2.17da}
&&F_{k,s,v}(y)
=(2y)^2\int_0^\iy\frac{t^{s-2m-1}e^{(1-v)t}}
 {e^t-(-1)^k}
 \frac{t^{2m}}{t^2+4y^2}dt\nonumber\\
&=&(2y)^2\int_0^\iy\frac{t^{s-2m-1}e^{(1-v)t}}
 {e^t-(-1)^k}
 \left(\sum_{j=1}^{m}
(-1)^{j-1}t^{2(m-j)}(2y)^{2j-2}
+(-1)^m\frac{(2y)^{2m}}
{t^2+4y^2}\right)dt.
\ena
Hence \eqref{E2.16} follows from
\eqref{E2.17da} and \eqref{E2.1}.

\subsection{Properties of $H$.}\label{S3.3}
Recall that $m_k$ is defined by
\eqref{E2.6d}.
To prove statement (d), we need certain
properties of the function
\ba
H(w)=H_{k,s,v,y}(w)
:=\frac{w^{s-1}e^{(1-v)w}}
{\left(e^w-(-1)^k\right)
\left(1+\left[w/(2y)\right]^2\right)},
\qquad y\in\R\setminus\{0\},
\quad v\in(0,1],
\ea
and its parts.
Here, we choose a branch of $w^{s-1}$
that is holomorphic on $\CC\setminus [0,\iy)$
and satisfies the condition
$ \lim_{b\to 0^+} (a+ib)^{s-1}=a^{s-1}$
for $a>0$.
We start with two simple properties.

\begin{property}\label{P2.8}
H is holomorphic in the complete angle
$0< \mathrm{arg}\, w< 2\pi$ with sides
$L_1=L_2=[0,\iy)$, except the points of the set
$E_{k,y}:=\{0\}\cup\{2iy,-2iy\}\cup
\{\pi i(2n+k):n\in\Z\}$,
which consists of the origin and
the simple poles of $H$.
\end{property}

\begin{property}\label{P2.9}
For $\vert w\vert>2\sqrt{2}\vert y\vert$,
\ba
h_1(w):=\left\vert w^{s-1}
/\left(1+\left[w/(2y)\right]^2\right)\right\vert
\le 2y^2 \vert w\vert^{\mathrm{Re}\,s-3},
\qquad s\in\CC.
\ea
\end{property}
\noindent
More relations  are discussed in the next three
properties.

\begin{property}\label{P2.10}
For $v\in(0,1]$ and
$\vert w\vert=\pi\left(2d+k-1\right),\,d\in\N$,
\beq\label{E2.22}
h_2(w):=\left\vert e^{(1-v)w}
\left(e^w-(-1)^k\right)^{-1}\right\vert
\le \left\{\begin{array}{ll}
2e^{-C_3(v)\,\left\vert\mathrm{Re}\,w
\right \vert},&v\in(0,1),\\
2, &v=1,
\end{array}\right.
\eeq
where
$C_3(v):=\left\{\begin{array}{ll}
v,&v\in(0,1/2),\\
1-v,&v\in [1/2,1].\end{array}\right.$
\end{property}
\proof
For $z=a+ib,\,a\in\R$, and $b\in\R$,
the following identities are valid:
\beq\label{E2.23}
\left\vert e^z\pm e^{-z}\right\vert^2
=\left(e^{a}\pm e^{-a}\right)^2\cos^2 b
+\left(e^{a}\mp e^{-a}\right)^2\sin^2 b
=e^{2a}+e^{-2a}\pm 2\cos 2b.
\eeq
First, it follows from \eqref{E2.23} that
\bna
&&\left\vert e^z-(-1)^k e^{-z}\right\vert
\ge\left(e^{2a}+e^{-2a}-2\right)^{1/2}
\ge (1/2)e^{\vert a\vert},\quad \vert a\vert\ge 1;
\label{E2.24}\\
&&\left\vert e^z-(-1)^k e^{-z}\right\vert
\ge e^{\vert a\vert}
\left\vert
\cos\left(b-\pi(k-1)/2\right)\right\vert,
\quad \vert a\vert< 1.\label{E2.25}
\ena
If $\vert z\vert=\pi\left(d+(k-1)/2\right),
\,d\in\N$, then for $\vert a\vert< 1$,
\beq\label{E2.26}
\left\vert\cos\left(b-\pi(k-1)/2
\right)\right\vert
=\left\vert\cos\left(
\vert z\vert-\sqrt{\vert z\vert^2-a^2}
\right)\right\vert
\ge \cos(1/\vert z\vert)
\ge \cos(2/\pi)
> 1/2.
\eeq
Next, combining \eqref{E2.24}, \eqref{E2.25},
and \eqref{E2.26}, we obtain
\beq\label{E2.27}
\left\vert e^z-(-1)^k e^{-z}\right\vert
\ge (1/2)e^{\vert\mathrm{Re}\,z \vert},\qquad
\vert z\vert=\pi\left(d-(1-k)/2\right),
\quad d\in\N.
\eeq
Finally, using \eqref{E2.27} for $z=w/2,\,
\vert w\vert=\pi\left(2d+k-1)\right),
\,d\in\N$, we have
\ba
h_2(w):=\frac{e^{(1/2-v)\mathrm{Re}\,w}}
{\left\vert e^{w/2}-(-1)^k e^{-w/2}\right\vert}
\le \frac{2}{e^{(v-1/2)\mathrm{Re}\,w
+\left\vert\mathrm{Re}\,w\right\vert/2}}.
\ea
Thus \eqref{E2.22} follows. \hfill $\Box$

\begin{property}\label{P2.11}
For $v\in(0,1]$ and $0<\vert w\vert
<\mathrm{min}\{2\vert y\vert,2\pi/3\}$,
\ba
\vert H(w)\vert
\le \frac{C_4(k)\vert w\vert^{
\mathrm{Re}\,s +k-2}e^{(1-v)\vert w\vert}}
{1-\left[\vert w\vert/(2y)\right]^2}.
\ea
\end{property}
\proof
It suffices to show that
\beq\label{E2.29}
\vert e^w-(-1)^k \vert \ge C_4^{-1}\vert w\vert^{1-k},
\qquad \vert w\vert\le 2\pi/3.
\eeq
For $w=a+ib,\,a\in\R,\,b\in\R$,
and $\vert w\vert\le 2\pi/3$,
the following relations are valid:
\bna\label{E2.30}
\vert e^w-(-1)^k \vert
&=&\left(\left(e^a-1\right)^2
+4e^{a}
\cos^2\left(b/2-\pi (k-1)/2\right)\right)^{1/2}
\nonumber\\
&\ge& \mathrm{max}\left\{\vert e^a-1\vert, 2e^{a/2}
\left\vert\cos\left(b/2+\pi (1-k)/2\right)
\right\vert\right\}.
\ena
It follows from \eqref{E2.30} that
\bna
&&\vert e^w-(-1)^k \vert
\ge 1-e^a
\ge 1/2
\ge C_5(k)\vert w\vert^{1-k},\qquad a<-\log 2;
\label{E2.31}\\
&&\vert e^w-(-1)^k \vert
\ge \mathrm{max}\left\{\vert a\vert/2,
(1/2)^{1/2}(2\vert b\vert/\pi)^{1-k}
\right\}\nonumber\\
&\ge& C_6(k)\left\{\begin{array}{ll}
1, &k=1,\\
\mathrm{max}\left\{\vert a\vert,\vert b\vert\right\},
&k=0,\end{array}\right.
\ge C_7(k)\vert w\vert^{1-k},
\qquad a\ge -\log 2,
\quad \vert b\vert\le 2\pi/3.\label{E2.32}
\ena
Thus \eqref{E2.29} follows from
\eqref{E2.31} and \eqref{E2.32}
with $C_4^{-1}:=\max\{C_5(k),C_7(k)\}$.
\hfill $\Box$

\begin{property}\label{P2.12}
Let $\mathrm{Re}\,s\in (1,3)$,
 if $k=0$ and $v\in(0,1)$;
$\mathrm{Re}\,s\in (1,2)$ if $k=0$ and $v=1$; and
$\mathrm{Re}\,s\in (0,2)$,
 if $k=1$ and $v\in(0,1]$.
Then
\bna
\lim_{\vep\to 0^+}\int_{K(\vep)}H(w)dw=0,\label{E2.33}\\
\lim_{d\to \iy}\int_{K(R_d)}H(w)dw=0,\label{E2.34}
\ena
where $R_d:=\pi\left(2d+k-1\right),\,d\in\N$.
\end{property}
\proof
Since the conditions of the property guarantee that
$\mathrm{Re}\,s+k-1>0$, we obtain from
Property \ref{P2.11}
\ba
\lim_{\vep\to 0^+}\int_{K(\vep)}\vert H(w)\vert dw
\le 2\pi C_4(k)\lim_{\vep\to 0^+}
\vep^{\mathrm{Re}\,s+k-1}=0.
\ea
Thus \eqref{E2.33} is established.
Next, for $\vert w\vert=R_d$, where $d\in\N$ is
large enough, we have from Properties
\ref{P2.9} and \ref{P2.10}
\beq\label{E2.35}
\vert H(w)\vert=h_1(w)h_2(w)
\le 4y^2\vert \vert w\vert^{\mathrm{Re}\,s-3}
\left\{\begin{array}{ll}
e^{-C_3(v)\vert w\cos (\mathrm{arg}\,w)\vert}, &v\in(0,1),\\
1, &v=1.
\end{array}\right.
\eeq
If $v=1$ and $\mathrm{Re}\,s\in(1,2)$ or
$\mathrm{Re}\,s\in(0,2)$, then by \eqref{E2.35},
\beq\label{E2.36}
\lim_{d\to \iy}\int_{K\left(R_d\right)}\vert H(w)\vert dw
\le 8\pi y^2\lim_{d\to\iy}R_d^{\mathrm{Re}\,s-2}=0.
\eeq
If $v\in(0,1)$ and $\mathrm{Re}\,s\in(1,3)$ or
$\mathrm{Re}\,s\in(0,2)$, then
using \eqref{E2.35} and Jordan's Lemma,
we obtain
\beq\label{E2.37}
\lim_{d\to \iy}\int_{K\left(R_d\right)}\vert H(w)\vert dw
\le 8\pi y^2\lim_{d\to\iy}R_d^{\mathrm{Re}\,s-3}=0.
\eeq
Thus \eqref{E2.34} follows from
\eqref{E2.36} and \eqref{E2.37} in all cases.
\hfill $\Box$

\subsection{Proof of Statement (d).}\label{S3.4}
We first prove  Lemma \ref{L2.7} (d)
under the conditions of Property \ref{P2.12} when
$m_k=0$, that is, we show that the following
relations hold true:
\bna\label{E2.38}
&&\sin(y+\pi k/2)F_{k,s,v}(y)
=f_{k,s,v}(y)-g_{k,s,v}(y)\nonumber\\
&&=\pi 2^{s-2}
\left(-\frac{\vert y\vert^s\sin[(2v-1)y-\pi k/2]}
{\sin(\pi s/2)}
-\frac{\vert y\vert^s\mathrm{sgn}\,y\cos[(2v-1)y-\pi k/2]}
{\cos(\pi s/2)}\right)\nonumber\\
&&+\frac{\pi 2^{s}\sin(y+\pi k/2)\,y^{2}}{\sin(\pi s)}
\sum_{n=1-k}^\iy
\frac{[\pi(n+k/2)]^{s-1}\sin[\pi(2n+k)v+\pi s/2]}
{y^2-[\pi(n+k/2)]^2}.
\ena
Indeed, using Properties \ref{P2.8} and \ref{P2.12}
and the Residue Theorem, we obtain
\ba
\left(1-e^{2\pi is}\right)\int_0^\iy H(w)dw
=\int_{L_1} H(w)dw
-e^{2\pi is}\int_{L_2} H(w)dw
=2\pi i\sum_{w\in E_{k,y}\setminus \{0\}}
\mathrm{Res}(H,w).
\ea
Hence
\beq\label{E2.39}
F_{k,s,v}(y)=\int_0^\iy H(w)dw
=\frac{2\pi i}{1-e^{2\pi is}}
\sum_{w\in E_{k,y}\setminus \{0\}}
\mathrm{Res}(H,w)
\eeq
(recall that $s\notin\N$).
Next, by straightforward calculations,
\bna\label{E2.40}
&&\frac{2\pi i}{1-e^{2\pi is}}
(\mathrm{Res}(H,2iy)+\mathrm{Res}(H,-2iy))
\nonumber\\
&&=\frac{\pi 2^{s-2}}
{\sin(y+\pi k/2)}
\left(-\frac{\vert y\vert^s\sin[(2v-1)y-\pi k/2]}
{\sin(\pi s/2)}
-\frac{\vert y\vert^s\mathrm{sgn}\,y\cos[(2v-1)y-\pi k/2]}
{\cos(\pi s/2)}\right)
\nonumber\\
&&=f_{k,s,v}(y)/\sin(y+\pi k/2)
\ena
and
\bna\label{E2.41}
&&\frac{2\pi i}{1-e^{2\pi is}}
\sum_{n\in\Z,n+2k\ne 0}
\mathrm{Res}(H,\pi i(2n+k))
\nonumber\\
&&=\frac{\pi 2^{s}i^s\,y^{2}}{1-e^{2\pi is}}
\sum_{n=1-k}^\iy
\frac{[\pi(n+k/2)]^{s-1}
\left(e^{-\pi i(2n+k)}
+(-1)^{s-1}e^{-\pi i(2n+k)}\right)}
{y^2-[\pi(n+k/2)]^2}\nonumber\\
&&=\frac{\pi 2^{s}\,y^{2}}{\sin(\pi s)}
\sum_{n=1-k}^\iy
\frac{[\pi(n+k/2)]^{s-1}\sin[\pi(2n+k)v+\pi s/2]}
{y^2-[\pi(n+k/2)]^2}
\nonumber\\
&&=-g_{k,s,v}(y)/\sin(y+\pi k/2).
\ena
Then it follows from
\eqref{E2.39}, \eqref{E2.40}, and \eqref{E2.41}
that the sum of series \eqref{E2.16a}
for $m_k=0$ is equal to
\ba
\frac{\sin(\pi s)}{\pi 2^{s}\,y^{2}}
\left(\sin(y+\pi k/2)\int_0^\iy H_{k,s,v,y}(w)dw
-f_{k,s,v}(y)\right),
\ea
that is, series \eqref{E2.16a} is convergent for
$y\in\R\setminus\{0\}$.
Note that the convergence of \eqref{E2.16a}
is trivial for $0<\mathrm{Re}\,s<2,\,s\ne 1,$
but it is not the case for $2<\mathrm{Re}\,s<3$.

Combining now \eqref{E2.39}, \eqref{E2.40},
and \eqref{E2.41},
we arrive at \eqref{E2.38} and \eqref{E2.17} as well
if the conditions of Property \ref{P2.12} are satisfied.

Furthermore, we extend \eqref{E2.17} to the case when
$k,\,s,$ and $v$ satisfy one of less restrictive conditions
\eqref{E2.6a}, \eqref{E2.6b}, and \eqref{E2.6c} of
Theorem \ref{T2.1}.
Since $k,\,s-2m_k,$ and $v$ satisfy the conditions
of Property \ref{P2.12},
we see from \eqref{E2.38} that
\bna\label{E2.42}
(-1)^{m_k}(2y)^{2m_k}F_{k,s-2m_k,v}(y)
&=&(-1)^{m_k}(2y)^{2m_k}\left(f_{k,s-2m_k,v}(y)
-g_{k,s-2m_k,v}(y)\right)\nonumber\\
&=&f_{k,s,v}(y)
-(-1)^{m_k}(2y)^{2m_k}g_{k,s-2m_k,v}(y).
\ena
Finally, using recurrence relation \eqref{E2.16}
with $m=m_k$, we arrive at \eqref{E2.17} from
\eqref{E2.42}.

Thus the proof of
Lemma \ref{L2.7} is completed.
\hfill $\Box$

Note that relations like \eqref{E2.17}
for $v=1/2$ and real $s$ were obtained in
\cite[Lemma 5 (a)]{G2002}.

\section{Properties of $g_{k,s,v}$.}\label{S4}
\setcounter{equation}{0}
Recall that $f_{k,s,v},\,
g_{k,s,v}$,
and $F_{k,s,v}$ are defined by
\eqref{E2.8}, \eqref{E2.7}, and \eqref{E2.14},
respectively.
In addition to Lemma \ref{L2.7},
 the proof of Theorem \ref{T2.1}
  is based
on certain properties of $g_{k,s,v}$.

\begin{lemma}\label{L2.13}
Let $k,\,s,$ and $v$ satisfy one of conditions
\eqref{E2.6a}, \eqref{E2.6b}, and \eqref{E2.6c} of
Theorem \ref{T2.1}.
Then the following statements are valid.\\
(a) $g_{k,s,v}$ is an EFET1
that satisfies the conditions
(C1) $f_{k,s,v}-g_{k,s,v}\in L_\iy(\R)$ and
(C2) $g_{k,s,v}^{(1-k)}(0)=0$.\\
(b) $g_{k,s,v}$ interpolates $f_{k,s,v}$
at the nodes $\{\pi(n+k/2)\}_{n\in\Z}$.
In addition, relations \eqref{E2.8a} hold true.\\
(c) $g_{k,s,v}$ is the only EFET1 that
interpolates $f_{k,s,v}$
at the nodes $\{\pi(n+k/2)\}_{n\in\Z}$
and satisfies
(C1)  and (C2).
\end{lemma}
\proof
(a) We first prove that $g_{k,s,v}\in B_1$.
It suffices to show that if
$k,\,s,$ and $v$ satisfy the conditions of
Property \ref{P2.12}, then $G(y):=g_{k,s,v}(y)/y^2$
defined by \eqref{E2.16a} for $m_k=0$ is an EFET1.
The traditional technique of proving such a result
 (see, e.g., \cite[Sect. 4.3]{T1963} and the
 proof of Corollary \ref{C2.10})
 can be applied here for $\mathrm{Re}\,s\in(0,2),
 \,s\ne 1$,
 but it is not applicable for
  $\mathrm{Re}\,s\in(2,3)$ since series
  \eqref{E2.16a} for $m_k=0$
  is not absolutely convergent. That is why
  we use a different approach based on the
  Fourier transform of $G$
  and on the generalized Paley-Wiener theorem.

  Note that by \eqref{E2.17},
  \eqref{E2.8}, and
  the boundedness of
 $F_{k,s,v}$ on $\R$,
 \beq\label{E2.49}
 \vert g_{k,s,v}(y)\vert
 \le \vert F_{k,s,v}(y)\vert
 +\vert f_{k,s,v}(y)\vert
 \le C_8(k,s,v)
 \left(1+\vert y\vert^{\mathrm{Re}\,s}\right).
 \eeq
 Next, if $\mathrm{Re}\,s\in (0,2),\,s\ne 1$ and
 $k=0$ or $k=1$, then series \eqref{E2.16a}
  for $m_k=0$ converges uniformly on $[-1,1]$.
  Therefore, by \eqref{E2.49},
\beq\label{E2.50}
\sup_{y\in\R}\vert G(y)\vert <\iy,\qquad
\mathrm{Re}\,s\in (0,2),\quad s\ne 1,\quad
k=0\,\, \mathrm{or}\, \,k=1.
\eeq
If $\mathrm{Re}\,s\in (2,3)$ and
 $k=0$, then by \eqref{E2.15a}
 and \eqref{E2.49},
 \beq\label{E2.51}
 \vert G(y)\vert
 \le C_9(k,s,v)
 \left\{\begin{array}{ll}
 \vert y\vert^{-2/\mathrm{Re}\,s},
 &\vert y\vert\le 1,\\
\vert y\vert^{\mathrm{Re}\,s-2},
 &\vert y\vert> 1,
 \end{array}\right.
  \qquad
\mathrm{Re}\,s\in (2,3),\quad
k=0.
\eeq
In addition, by Lemma \ref{L2.7} (d), the function
$G$ is the limit of continuous functions
on $\R$. Hence $G$ is a measurable
function on $\R$, and it is locally
integrable on $\R$ by estimates
\eqref{E2.50} and \eqref{E2.51}.
Thus $G$ generates the tempered
 distribution $G$ by the formula
 $ \langle G,\psi \rangle:=\int_\R G(y)\psi(y)dy$
 for every test function $\psi$ from the
 Schwartz class $S(\R)$.

 Its distributional Fourier transform is given by the formula
 \beq\label{E2.52}
 \widehat{G}(u)
 =\frac{\pi 2^{s}}{\sin(\pi s)}
 \sum_{n=1-k}^\iy
[\pi(n+k/2)]^{s-1}\sin[\pi(2n+k)v+\pi s/2]\,
\widehat{h}_{k,n}(u),
\eeq
where
\beq\label{E2.52a}
h_{k,n}(y):=\frac{\sin(y+\pi k/2)}
{y^2-[\pi(n+k/2)]^2}
 \eeq
and
\beq\label{E2.53}
\widehat{h}_{k,n}(u)
=\left\{\begin{array}{ll}
\frac{(-1)^{n+k}\,i^{1-k}}{n+k/2}
\sin[\pi(n+k/2)u+\pi k/2], &\vert u\vert\le 1,\\
0,&\vert u\vert> 1.
\end{array}\right.
\eeq
To prove \eqref{E2.53}, we first note that
the following formulae for the Fourier transform:
\beq\label{E2.54}
\int_{\R}
\frac{\sin(y+\pi k/2)}
{y^2+b^2}\,e^{iuy}dy
=\frac{\pi i^{1-k}}{b}
\left\{\begin{array}{ll}
e^{-b}\sinh bu, &k=0,\quad \vert u\vert\le 1,\\
e^{-b\vert u\vert}\mathrm{sgn}  \,u\sinh b, &k=0,
\quad \vert u\vert> 1,\\
e^{-b}\cosh bu, &k=1,\quad \vert u\vert\le 1,\\
e^{-b\vert u\vert}\cosh b, &k=1,\quad \vert u\vert> 1,
\end{array}\right.\quad \mathrm{Re}\,b>0,
\eeq
immediately follow from the Laplace-type integral
\ba
\int_{\R}\frac{\cos ay}{y^2+b^2}\,dy=\frac{\pi}{b}\,
e^{-\vert a\vert b},\qquad \mathrm{Re}\,b>0,
\quad a\in\R,
\ea
(see, e.g., \cite[Eqn. 3.723.2]{GR1980}).
Next, setting $b=b_{k,n}(\de):=\de+i\pi (n+k/2)$
 in \eqref{E2.54},
we see that for all $y\in\R,\,\de\in(0,1)$, and
$\vert n\vert\ge 1-k$,
\bna\label{E2.54a}
\left\vert \frac{\sin(y+\pi k/2)}
{y^2+b_{k,n}^2(\de)}\right\vert
&\le& \frac{\left\vert \sin(y+\pi k/2)
\right\vert}
{\max\left\{\left\vert y^2-[\pi(n+k/2)]^2
+\de^2\right\vert,
2\pi\vert n+k/2\vert\de\right\}}\nonumber\\
&\le& \frac{1}{\min\left\{1/2,\pi
\vert n+k/2\vert\right\}}
\left\vert
\frac{\sin(y+\pi k/2)}
{y^2-[\pi(n+k/2)]^2}\right\vert
\le\frac{4}{\pi}.
\ena
Note that the second inequality of \eqref{E2.54a}
is proved by the following two cases:
\ba
\de\in
\left(0,\left\vert y^2
-[\pi(n+k/2)]^2\right\vert/2\right]\quad
\mathrm{and}\quad \de\in
\left(\left\vert y^2-[\pi(n+k/2)]^2\right\vert/2,1
\right).
\ea
Hence by the Dominated Convergence Theorem,
\beq\label{E2.55}
\widehat{h}_{k,n}(u)
=\lim_{\de\to 0^+}
\int_{\R}
\frac{\sin(y+\pi k/2)}
{y^2+b_{k,n}^2(\de)}\,e^{iuy}dy,
\eeq
and \eqref{E2.53} follows from \eqref{E2.54}
 and \eqref{E2.55}.

 Then combining \eqref{E2.52} and \eqref{E2.53},
 we see that the corresponding trigonometric series
 converges to a tempered distribution on $[-1,1]$
 and the support of the tempered distribution
 $\widehat{G}$ is a subset of $[-1,1]$.
 Finally using the generalized Paley-Wiener
 theorem
 (see, e.g., \cite[Theorem 7.2.3]{S2003}),
 we arrive at $G\in B_1$ and $g_{k,s,v}\in B_1$.

 Finally, the condition
 $g_{k,s,v}^{(1-k)}(0)=0$ is a consequence of
 representation \eqref{E2.7}, and
 it immediately follows from
 \eqref{E2.17} that
 $f_{k,s,v}-g_{k,s,v}\in L_\iy(\R)$.
 \vspace{.12in}\\
 (b) $g_{k,s,v}$ interpolates $f_{k,s,v}$
at the nodes $\{\pi(n+k/2)\}_{\vert n\vert=1-k}^\iy$
by  \eqref{E2.17}.
In addition, $f_{0,s,v}(0)=g_{0,s,v}(0)=0$,
by \eqref{E2.8} and \eqref{E2.7}.
To verify equalities \eqref{E2.8a}, we first find
$f_{k,s,v}[\pi(n+k/2)]$ for
$\vert n\vert\ge 1-k$ by a simple calculation and then
use \eqref{E2.17} again.
Note that \eqref{E2.8a} is valid for $n=k=0$ as well.
It is also possible to prove \eqref{E2.8a} without
using \eqref{E2.17} (at least for $m_k=0$) by
a straightforward calculation of $g_{k,s,v}[\pi(n+k/2)]$
that coincides with $f_{k,s,v}[\pi(n+k/2)]$
for $n\in\Z$.
 \vspace{.12in}\\
 (c) Assume that an EFET1  $g$
 interpolates $f_{k,s,v}$
at the nodes $\{\pi(n+k/2)\}_{n\in\Z}$,
and, in addition,
$f_{k,s,v}-g\in L_\iy(\R)$ and
$g^{(1-k)}(0)=0$.
Next, denoting $\GG:=g_{k,s,v}-g$, we obtain
from statements (a) and (b) of Lemma \ref{L2.13} that
$\GG\in B_1\cap L_\iy(\R)$ and $\GG(\pi(n+k/2))=0$
for $n\in\Z.$
It is well known
(see, e.g., \cite[Sect. 4.3.1]{T1963})
that there exists a constant $C$ such that
$\GG(y)=C\sin(y+\pi k/2)$.
Since $0=\GG^{(1-k)}(0)=C$, we arrive at $g=g_{k,s,v}$.
This completes the proof of Lemma \ref{L2.13}.
\hfill $\Box$

\section{Proofs of Main Results}\label{S5}
\setcounter{equation}{0}
Here, we prove Theorems \ref{T2.1} and \ref{T2.7}
 and Corollaries
\ref{C2.2}, \ref{C2.8}, \ref{C2.10},
 and \ref{C2.9}.
\vspace{.12in}\\
\emph{Proof of Theorem \ref{T2.1}.}
Statement (a) of the theorem
immediately follows from
 Lemma \ref{L2.13},
  while statement (c)
 is the direct consequence of relations
  \eqref{E2.15} and \eqref{E2.17} of
  Lemma \ref{L2.7}.
  Equalities \eqref{E2.8a} of statement (b)
   are proved in Lemma \ref{L2.13} (b).
  \hfill $\Box$\vspace{.12in}\\
\emph{Proof of Corollary \ref{C2.2}.}
(a) If $\Phi\left((-1)^k,s,v\right) =0$,
then by \eqref{E2.15} and \eqref{E2.17}
of Lemma \ref{L2.7},
\beq\label{E2.56}
\left\vert f_{k,s,v}(y)-g_{k,s,v}(y)
\right\vert
\le \left\vert F_{k,s,v}(y)\right\vert
\le C_1(k,s,v)y^{-2},
\qquad y\ne 0.
\eeq
Thus statement (a) follows from
\eqref{E2.56}.\\
(b) The statement is proved by contradiction.
Assume that there exists $p\in (1/2,\iy)$ such that
\beq\label{E2.57}
\left\|f_{k,s,v}-g_{k,s,v}\right\|_{L_p(\R)}<\iy,
\eeq
and, in addition,
$\Phi\left((-1)^k,s,v\right) \ne 0$.
For $N>0$ let us define a set
$M_N:=\{y\in [-N,N]:
\vert \sin(y+\pi k/2)\vert\ge 1/2\}$.
Then by statement (c) of Theorem \ref{T2.1},
there exists $y_0>0$ such that for all
$N> y_0$,
\ba
\inf_{y\in M_N\setminus
[- y_0, y_0]}
\left\vert f_{k,s,v}(y)-g_{k,s,v}(y)
\right\vert
\ge
\left\vert \G(s)\Phi
\left((-1)^k,s,v\right)\right\vert/4.
\ea
Hence
\ba
\left\|f_{k,s,v}-g_{k,s,v}\right\|_{L_p(\R)}
&\ge& \limsup_{N\to\iy}
\left\|f_{k,s,v}-g_{k,s,v}\right\|_
{L_p\left(M_N\setminus
[-y_0,y_0]\right)}\\
&\ge& \left(\left\vert\G(s)
\Phi\left((-1)^k,s,v\right)\right\vert/4\right)
\lim_{N\to\iy}\left\vert M_N\setminus
[-y_0,y_0]\right\vert
^{1/p}
=\iy,
\ea
which is in contradiction to \eqref{E2.57}.
\hfill $\Box$\vspace{.12in}\\
\emph{Proof of Theorem \ref{T2.7}.}
The proof is similar to those of
Lemmas \ref{L2.7} and \ref{L2.13}.
We first assume  that $\mathrm{Re}\,s\in
(0,2),\,s\ne 1$, i.e., $m_1=0$. Then it follows from
\eqref{E2.13a} that
\bna\label{E2.58}
\g_{s,q}(y)
&=&-\frac{\pi 2^{s}y^{2}}{\sin(\pi s)}
\sin y
\sum_{n\in \N\setminus E(\chi)}
\frac{(\pi n)^{s-1}
\sum_{l=1}^{q-1}\chi(l)\sin[2\pi n l/q+\pi s/2]}
{y^2-(\pi n)^2}\nonumber\\
&=&-\frac{\pi 2^{s}y^{2}}{\sin(\pi s)}
\sin y
\sum_{n=1}^\iy
\frac{(\pi n)^{s-1}
\sum_{l=1}^{q-1}\chi(l)\sin[2\pi n l/q+\pi s/2]}
{y^2-(\pi n)^2},
\ena
since
$\sum_{l=1}^{q-1}\chi(l)
\sin[2\pi n l/q+\pi s/2]=0,\,n\in E(\chi)$, by
\eqref{E2.6.2} and
the definition of $E(\chi)$.
Next, the series, representing $g_{0,s,l/q},\,
l=1,\ldots, q-1,$ in \eqref{E2.7},
 as well as the series,
representing $\g_{s,q}$ in \eqref{E2.58},
absolutely converge for $\mathrm{Re}\,s\in
(0,2),\,s\ne 1$.
Therefore, the following representations
for $\g_{s,q}$ and $\vphi_{s,q}$ are valid:
\beq\label{E2.59}
\g_{s,q}(y)=\sum_{l=1}^{q-1}\chi(l)g_{0,s,l/q},
\quad
\vphi_{s,q}(y)=\sum_{l=1}^{q-1}\chi(l)f_{0,s,l/q}
,\qquad \mathrm{Re}\,s\in (0,2),\quad s\ne 1.
\eeq
Despite equalities \eqref{E2.59}
Theorem \ref{T2.7} for $\mathrm{Re}\,s\in(0,2)$
 does not follow directly
from Theorem \ref{T2.1} for $k=0$ and $m_0=0$
 because $\mathrm{Re}\,s\in(1,3)$
in this case of Theorem \ref{T2.1}.
That is why we briefly discuss
below all major steps of
Lemmas \ref{L2.7} and \ref{L2.13}
that are used in the proof of  Theorem \ref{T2.7}.

To prove Theorem \ref{T2.7} (a),
we study properties of
the integral
\beq\label{E2.60}
F_{s,q}^*(y):=
\int_0^\iy\frac{t^{s-1}
\sum_{l=1}^{q-1}\chi(l)e^{(1-l/q)t}}
 {\left(e^t-1\right)
 \left(1+[t/(2y)]^2\right)}dt
 =\sum_{l=1}^{q-1}\chi(l)F_{0,s,l/q}(y),
\eeq
where $y\in\R\setminus\{0\},\,
 \mathrm{Re}\,s\in (0,2),\,s\ne 1$,
 and $F_{0,s,l/q}$
 is defined by \eqref{E2.14}.
 Then it is easy to see from \eqref{E2.59}
 and \eqref{E2.6.5}
 (cf. the proof of Lemma
 \ref{L2.7}) (a) that
 \beq\label{E2.61}
 \left\vert F_{s,q}^*(y)-\G(s)q^s L(s,\chi)
 \right\vert
 \le  C_{10}(s,q)y^{-2}.
 \eeq
 Next, we show that the following equality holds true:
 \beq\label{E2.62}
\sin y\,F_{s,q}^*(y)
=\vphi_{s,q}(y)-\g_{s,q}(y).
\eeq
The proof is similar to the one of
Lemma \ref{L2.7} (d).
Let us define the function
\ba
H^*(w)=H_{s,q,y}^*(w)
:=\frac{w^{s-1}\sum_{l=1}^{q-1}
\chi(l)e^{(1-l/q)w}}
{\left(e^w-1\right)
\left(1+\left[w/(2y)\right]^2\right)},
\qquad y\in\R\setminus\{0\},
\ea
that is
holomorphic in the complete angle
$0< \mathrm{arg}\, w< 2\pi$,
except the points of the set
$E_{y}^*:=\{0\}\cup\{2iy,-2iy\}\cup
\{2\pi i n:n\in\Z\}$,
which consists of the origin and
the simple poles of $H^*$.

We also need the following estimates of
$H^*(w)$ for $y\in\R\setminus\{0\}$:
\bna
&&\vert H^*(w)\vert
\le C_{11} y^2\vert \vert w\vert^{\mathrm{Re}\,s-3}
e^{-C_{12} \vert w\cos (\mathrm{arg}\,w)\vert},
\quad \vert w\vert=\pi(2d-1),
\quad d\in\N \label{E2.63};\\
&&\vert H^*(w)\vert
\le \frac{C_{13}\vert w\vert^{
\mathrm{Re}\,s -1}}
{1-\left[\vert w\vert/(2y)\right]^2},
\qquad
0<\vert w\vert
<\mathrm{min}\{2\vert y\vert,2\pi/3\},
\label{E2.64}
\ena
where the constants $C_{11}, C_{12}$,
and $C_{13}$ are independent
of $w$.
Inequality \eqref{E2.63} follows from
Properties  \ref{P2.9} and \ref{P2.10},
while \eqref{E2.64} is a consequence of
 inequality \eqref{E2.29} for $k=0$
 and  elementary estimates
\ba
\left\vert
\sum_{l=1}^{q-1}\chi(l)e^{(1-l/q)w}\right\vert
\le \sum_{l=1}^{q-1}
\left\vert e^{(1-l/q)w}-1\right\vert
\le (q-1)\,e^{2\pi/3}\,\vert w\vert,\qquad
\vert w\vert\le 2\pi/3.
\ea
Then using \eqref{E2.63}, \eqref{E2.64},
and Jordan's Lemma
(similarly to the proof of Property
\ref{P2.12}), we arrive at
\beq\label{E2.65}
\lim_{\vep\to 0^+}\int_{K(\vep)}H^*(w)dw
=\lim_{d\to \iy}\int_{K(\pi(2d-1))}H^*(w)dw=0.
\eeq
Furthermore, using \eqref{E2.65}
and the Residue Theorem
(similarly to the proof of Lemma
\ref{L2.7} (d)) we obtain
\ba
F_{s,q}^*(y)=\int_0^\iy H^*(w)dw
=\frac{2\pi i}{1-e^{2\pi is}}
\sum_{w\in E_{y}^*\setminus \{0\}}
\mathrm{Res}(H^*,w)
=\frac
{\vphi_{s,q}(y)-\g_{s,q}(y)}
{\sin y}.
\ea
Thus equality \eqref{E2.62} is established.

Finally, we discuss certain properties of
$\g_{s,q}$.  If $\mathrm{Re}\,s\in(0,2)$,
then $\g_{s,q}$ given by \eqref{E2.58}
belongs to $B_1$.
To prove this statement, it suffices to show that
$G^*(y):=\g_{s,q}(y)/y^2$
is an EFET1.
Indeed,
we first note that by \eqref{E2.62},
  \eqref{E2.13b}, and
  the boundedness of
 $F_{s,q}^*$ on $\R$,
 \beq\label{E2.66}
 \vert \g_{s,q}(y)\vert
 \le \vert F_{s,q}^* (y)\vert
 +\vert \vphi_{s,q}(y)\vert
 \le C_{14}(s,q)
 \left(1+\vert y\vert^{\mathrm{Re}\,s}\right).
 \eeq
 Next, since $\mathrm{Re}\,s\in (0,2)$,
 series \eqref{E2.58}
  converges uniformly on $[-1,1]$.
  Therefore by \eqref{E2.66},
\beq\label{E2.67}
\sup_{y\in\R}\vert G^*(y)\vert <\iy,\qquad
\mathrm{Re}\,s\in (0,2).
\eeq
In addition, the function
$\g_{s,q}(y)$ is the limit of continuous functions
on $\R$.
Hence $G^*$ is a measurable
function on $\R$ and inequality \eqref{E2.67}
holds true.
Thus $G^*$ is locally
 integrable on $\R$, and it generates the tempered
 distribution $G^*$ by the formula
 $(G^*,\psi):=\int_\R G^*(y)\psi(y)dy$
 for every test function $\psi$ from the
 Schwartz class $S(\R)$.

 Its distributional Fourier transform is given
  by the formula
 \beq\label{E2.68}
 \widehat{G}^*(u)
 =\frac{\pi 2^{s}}{\sin(\pi s)}
 \sum_{n=1}^\iy
[\pi n]^{s-1}
\sum_{l=1}^{q-1}\chi(l)\sin[2\pi n l/q+\pi s/2]\,
\widehat{h}_{0,n}(u),
\eeq
where $h_{0,n}$ and $\Hat{h}_{0,n}$ are given in
\eqref{E2.52a} and \eqref{E2.53}, respectively.

Then combining \eqref{E2.68} and \eqref{E2.53},
 we see that the corresponding trigonometric series
 converges to a tempered distribution on $[-1,1]$
 and the support of the tempered distribution
 $\widehat{G}^*$ is a subset of $[-1,1]$.
 Finally, using the generalized Paley-Wiener
 theorem
 (see, e.g., \cite[Theorem 7.2.3]{S2003}),
 we arrive at $G^*\in B_1$ and $\g_{s,q}\in B_1$.

 In addition, it is easy to verify equalities
 \eqref{E2.13c} (see \eqref{E2.8a}) and to show
 that conditions
 (C1*) and (C2*) are satisfied. The uniqueness of
 $\g_{s,q}$  that interpolates
$\vphi_{s,q}$
at the nodes $\{\pi n\}_{n\in\Z}$
and satisfies conditions (C1*) and (C2*) can be proved
similarly to the proof of Lemma \ref{L2.13} (c).

Next, relation \eqref{E2.13d} follows from
\eqref{E2.61} and \eqref{E2.62}. Thus Theorem \ref{T2.7}
is established for $\mathrm{Re}\,s\in(0,2)$.

Let now
$\mathrm{Re}\,s\in (0,\iy),\,
 \mathrm{Re}\,s\ne 2,\,4,\ldots $, and
 $s\notin\N.$
 Recall that
 $m_1=\lfloor (\mbox{Re}\,s)/2\rfloor$.
Since $s-2m_1\in(0,2)$,
we can replace $s$ by $s-2m_1$
in \eqref{E2.62} and obtain the identity
\bna\label{E2.69}
(-1)^{m_1}(2y)^{2m_1}F^*_{s-2m_1,q}(y)
&=&(-1)^{m_1}(2y)^{2m_1}\left(\vphi_{s-2m_1,q}(y)
-\g_{s-2m_1,q}(y)\right)\nonumber\\
&=&\vphi_{s,q}(y)
-(-1)^{m_1}(2y)^{2m_1}\g_{s-2m_1,q}(y).
\ena
Finally, using \eqref{E2.60} and
recurrence relation \eqref{E2.16}
with $m=m_1$, we arrive at \eqref{E2.13d} from
\eqref{E2.69}.
Thus the proof of
Theorem \ref{T2.7} is completed.
\hfill $\Box$\vspace{.12in}\\
\emph{Proof of Corollary \ref{C2.8}.}
The proof of the corollary is based on
inequality \eqref{E2.61} and statement (c) of
Theorem \ref{T2.7} similarly to the proof of
Corollary \ref{C2.2}.
\hfill $\Box$\vspace{.12in}\\
\emph{Proof of Corollary \ref{C2.10}.}
Note that by Proposition \ref{P1.1} (a),
$E(\chi(\cdot,3))=\{3d\pm 1\}_{d\in\Z}$.
Then formulae \eqref{E2.14a}, \eqref{E2.14b}, and
\eqref{E2.14d} immediately follow from the relations
\beq\label{E2.70}
\g_{s,3}^*(y)=\g_{s,3}(3y/2)/\sin(y/2),\qquad
\vphi_{s,3}^*(y)=\vphi_{s,3}(3y/2)/\sin(y/2),
\eeq
and equalities \eqref{E2.13a}, \eqref{E2.13b}, and
\eqref{E2.13d} for $q=3$, respectively.

It remains to prove \eqref{E2.14c}
and the statement that $\g_{s,3}^*$, given by
\eqref{E2.14a},
is the only EFET1 that interpolates
$\vphi_{s,3}^*$
at the nodes $\{2\pi (d\pm 1/3)\}_{d\in\Z}$
and satisfies conditions
(C1*) and (C2*).

We first prove that $\g_{s,3}^*\in B_1$.
We use the traditional technique
(see, e.g., \cite[Sect. 4.3]{T1963}) for the proof
because we do not know as to whether
the Fourier method used in the proof of
Lemma \ref{L2.13} (a) can be applied in this case.

It follows from \eqref{E2.70}
that $\g_{s,3}^*$ is an entire function
since $\g_{s,3}$ is entire
 by Theorem \ref{T2.7} (a)
 and all zeros of $\sin(y/2)$ are zeros of
 $\g_{s,3}(3y/2)$.
Therefore, to prove that $\g_{s,3}^*\in B_1$, it
 suffices to estimate the function series
\beq\label{E2.71}
\psi_r(z):=\sum_{d=0}^\iy (d+r/3)^{\mathrm{Re}\,s-1}
\left\vert\frac{1+2\cos z}
{z^2-[2\pi(d+r/3)]^2}\right\vert,
\eeq
where $ r=1$ or $r=2,\, z=x+iy\in\CC,$ and
$\mathrm{Re}\,s\in(0,2)$.
Since $\psi_r$ is an even function, we can assume,
without loss of generality, that
$x=\mathrm{Re}\,z\ge 0$.
Next, for every $x\ge 0$
there exists $d_x\in\Z_+$ such that
$\left\vert x-2\pi(d_x+r/3)\right\vert\le 4\pi/3$.
Then
\bna\label{E2.72}
\psi_r(z)
&=&(d_x+r/3)^{\mathrm{Re}\,s-1}
\left\vert\frac{1+2\cos z}
{z^2-[2\pi(d_x+r/3)]^2}\right\vert\nonumber\\
&+&\sum_{d=0,\,d\ne d_x}^\iy (d+r/3)^{\mathrm{Re}\,s-1}
\left\vert\frac{1+2\cos z}
{z^2-[2\pi(d+r/3)]^2}\right\vert
=I_1(z)+I_2(z),
\ena
where
\beq\label{E2.73}
I_1(z)
\le 3\pi^{-2}\left(d_x+r/3\right)^{\mathrm{Re}\,s-1}
e^{\vert y\vert}
\le \left(\vert z\vert+1\right)^{\mathrm{Re}\,s}
e^{\vert \mathrm{Im\,z}\vert},\qquad
\vert y\vert=\vert\mathrm{Im}\,z\vert\ge\pi.
\eeq
Furthermore, if $\vert y\vert<\pi$, then
setting $w:=z-2\pi(d_x+r/3)$, we see that
$\vert w\vert< (5/3)\pi$.
Therefore,
\bna\label{E2.74}
\left\vert\frac{1+2\cos z}
{z-2\pi(d_x+r/3)}\right\vert
&=&4\left\vert\frac{\sin(w/2)}{w}\right\vert
\left\vert\cos(w/2+(-1)^{r+1}\pi/6)\right\vert
\nonumber\\
&\le& \frac{4\sinh(\vert w\vert/2)}
{\vert w\vert}
\cosh(\vert w\vert/2+\pi/6)
<C_{15}, \qquad \vert y\vert <\pi,
\ena
where $C_{15}\in (1,61)$ is an absolute constant.
Since
$\vert z+2\pi(d_x+r/3)\vert^{-1}\le
[2\pi(d_x+r/3)]^{-1}$ for $x\ge 0$,
we obtain from \eqref{E2.74}
\beq\label{E2.75}
I_1(z)\le C_{15}(d_x+r/3)^{\mathrm{Re}\,s-2}
<C_{15},\qquad \vert y\vert<\pi.
\eeq
Combining \eqref{E2.73} and \eqref{E2.75},
we arrive at the estimate
\beq\label{E2.76}
I_1(z)\le C_{15}
\left(\vert z\vert+1\right)^{\mathrm{Re}\,s}
e^{\vert \mathrm{Im\,z}\vert},\qquad
z\in\CC.
\eeq
Next,
\bna\label{E2.77}
I_2(z)&\le&
3e^{\vert y\vert}
\sum_{d=0,\,d\ne d_x}^\iy (d+r/3)^{\mathrm{Re}\,s-1}
\vert x^2-[2\pi(d+r/3)]^2\vert^{-1}\nonumber\\
&\le& 3(2\pi)^{-2}e^{\vert y\vert}
\sum_{d=0,\,d\ne d_x}^\iy (d+r/3)^{\mathrm{Re}\,s-1}
\left(\vert d-d_x \vert-2/3\right)^{-1}
\left(\vert d+d_x \vert+2(r-1)/3\right)^{-1}\nonumber\\
&\le& 9(2\pi)^{-2}e^{\vert y\vert}
\sum_{d=0,\,d\ne d_x}^\iy (d+r/3)^{\mathrm{Re}\,s-1}
\left\vert d^2-d_x^2 \right\vert^{-1}\nonumber\\
&\le& C_{16}(s)e^{\vert y\vert}\left(
\sum_{d=0}^{d_x-1} (d+1)^{\mathrm{Re}\,s-1}
\left(d_x^2-d^2 \right)^{-1}
+\sum_{d=d_x+1}^\iy (d+1)^{\mathrm{Re}\,s-1}
\left(d^2-d_x^2 \right)^{-1}
\right)\nonumber\\
&=& C_{16}(s)e^{\vert y\vert}\left(J_{1,d_x}
+J_{2,d_x}\right),
\ena
where
$C_{16}(s):=9(2\pi)^{-2}
\max \{1,3^{1-\mathrm{Re}\,s}\}$ and
$J_{1,0}=0$. Then
\bna\label{E2.78}
J_{1,d_x}&=&\sum_{\nu=1}^{d_x}
\frac{\left(d_x-\nu+1\right)^{\mathrm{Re}\,s-1}}
{\nu\left(2d_x-\nu\right)}
\le \left\{\begin{array}{ll}
d_x^{-1}
\sum_{\nu=1}^{d_x}1/\nu, &\mathrm{Re}\,s\in(0,1),\\
d_x^{\mathrm{Re}\,s-2}
\sum_{\nu=1}^{d_x}1/\nu, &\mathrm{Re}\,s\in [1,2),
\end{array}\right.\nonumber\\
&\le& \left\{\begin{array}{ll}
d_x^{-1}
\left(\log d_x+1\right), &\mathrm{Re}\,s\in(0,1),\\
d_x^{\mathrm{Re}\,s-2}
\left(\log d_x+1\right), &\mathrm{Re}\,s\in [1,2),
\end{array}\right.
\le C_{17}(s),\quad d_x\in\N,
\ena
and
\beq\label{E2.79}
J_{2,0}=\sum_{d=1}^{\iy}
\frac{\left(d+1\right)^{\mathrm{Re}\,s-1}}
{d^{2}};\qquad
J_{2,d_x}=\sum_{\nu=1}^{\iy}
\frac{\left(d_x+\nu+1\right)^{\mathrm{Re}\,s-1}}
{\nu\left(2d_x+\nu\right)}
\le \sum_{\nu=1}^{\iy}\nu^{\mathrm{Re}\,s-3},\quad
d_x\in\N.
\eeq
Combining \eqref{E2.77}, \eqref{E2.78},
and \eqref{E2.79},
we arrive at the estimate
\beq\label{E2.80}
I_2(z)\le C_{18}(s)
e^{\vert \mathrm{Im\,z}\vert},
\qquad z\in\CC.
\eeq
Finally, taking account of \eqref{E2.72},
 \eqref{E2.76}, and \eqref{E2.80}, we obtain
 the estimate
 \ba
\psi_r(z)\le
C\left(\vert z\vert+1\right)^{\mathrm{Re}\,s}
e^{\vert \mathrm{Im\,z}\vert},
\qquad z\in\CC,
\ea
where $C$ is independent of $z$.
Hence \eqref{E2.71} and \eqref{E2.14a}
 show that $g_{s,3}^*\in B_1$.

 Next, it immediately follows from
 representation \eqref{E2.14a} that
 the condition
 $\g_{s,3}^*(0)=0$ is satisfied.
 In addition,
 $\vphi_{s,3}^*-\g_{s,3}^*\in L_\iy(\R)$.
 Indeed, recall that by \eqref{E2.69},
 equality \eqref{E2.62} is valid for
 $\mathrm{Re}\,s\in (0,\iy),\,
 \mathrm{Re}\,s\ne 2,\,4,\ldots $, and
 $s\notin\N.$
 Then it follows from \eqref{E2.62} and
 \eqref{E2.70} that
 \beq\label{E2.81}
 (1+2\cos y)\left(
 F_{0,s,1/3}(3y/2)-F_{0,s,2/3}(3y/2)\right)
 =\vphi_{s,3}^*(y)-\g_{s,3}^*(y).
 \eeq
 Hence condition (C1*) is satisfied.
 In addition, \eqref{E2.14c} is a consequence of
 \eqref{E2.14b} and \eqref{E2.81}.
 It is also possible to prove \eqref{E2.14c} without
using \eqref{E2.81} (at least for $m_1=0$) by
a straightforward calculation of
$\g_{s,3}[2\pi (d\pm 1/3)]$
that coincides with $\vphi_{s,3}[2\pi (d\pm 1/3)]$
for $d\in\Z$.

 Finally, let us assume that an EFET1  $g$
 interpolates $\vphi_{s,3}^*$
at the nodes $\{2\pi (d\pm 1/3)\}_{d\in\Z}$
and, in addition,
$\vphi_{s,3}^*-g\in L_\iy(\R)$ and
$g(0)=0$.
Next, denoting $\GG:=\g_{s,3}^*-g$, we obtain
 that
$\GG\in B_1\cap L_\iy(\R)$ and $\GG(0)=0$.
In addition,
$\GG(2\pi (d\pm 1/3))=0,\,{d\in\Z}$.
Then the function
$\GG_1(y):=\GG(y)\sin(y/2)$ belongs to
$B_{3/2}\cap L_\iy(\R)$ and
$\GG_1\left(\frac{\pi n}{3/2}\right)=0,\,{n\in\Z}$.
It is well known
(see, e.g., \cite[Sect. 4.3.1]{T1963})
that there exists a constant $C$ such that
$\GG_1(y)=C\sin(3y/2)$, that is,
$\GG(y)=C(1+2\cos y)$.
Since $0=\GG(0)=3C$, we arrive at $g=\g_{s,3}^*$.
This completes the proof of Corollary \ref{C2.10}.
 \hfill $\Box$\vspace{.12in}\\
\emph{Proof of Corollary \ref{C2.9}.}
Formulae \eqref{E2.14e} immediately follow
 from \eqref{E2.13a} and \eqref{E2.13b},
 while \eqref{E2.14f}
 is a consequence of \eqref{E2.6.3b} and
  \eqref{E2.13}.
  In addition, \eqref{E2.14f} follows as well from
  \eqref{E2.13d} for $q=4$.
 \hfill $\Box$

\end{document}